\documentclass{amsart}

\usepackage{amscd,amssymb}
\usepackage[all]{xy}

\usepackage{hyperref}

\def\Z{{\mathbb Z}}
\def\Q{{\mathbb Q}}
\def\R{{\mathbb R}}
\def\C{{\mathbb C}}
\def\F{{\mathbb F}}
\def\P{{\mathbb P}}
\def\U{{\mathbb U}}
\def\V{{\mathbb V}}
\def\disk{{\mathbb D}}

\def\A{{\mathcal A}}
\def\B{{\mathcal B}}
\def\E{{\mathcal E}}

\def\H{{\mathcal H}}
\def\J{{\mathcal J}}
\def\L{{\mathcal L}}
\def\M{{\mathcal M}}
\def\N{{\mathcal N}}
\def\O{{\mathcal O}}
\def\cC{{\mathcal C}}
\def\h{{\mathfrak h}}

\def\b{{\beta}}
\def\d{{\delta}}
\def\s{{\sigma}}
\def\w{{\omega}}

\def\G{{\Gamma}}
\def\D{{\Delta}}

\def\Mtilde{\widetilde{\M}}

\def\Bbar{\overline{\B}}
\def\Fbar{\overline{F}}
\def\Hbar{\overline{\H}}
\def\Mbar{\overline{\M}}
\def\Wbar{\overline{W}}

\def\Bhat{{\widehat{\B}}}

\def\Jdual{\check{\J}}
\def\Udual{\check{U}}
\def\bUdual{\check{\U}}
\def\Vdual{\check{V}}
\def\bVdual{\check{\V}}

\def\qtilde{\w_q}
\def\utilde{\w_u}
\def\wtilde{\tilde\w}

\def\ebar{\overline{e}}
\def\fbar{\overline{f}}
\def\gbar{\overline{g}}
\def\zbar{\overline{z}}
\def\ubar{\overline{u}}
\def\vbar{\overline{v}}

\def\ep{{\epsilon}}
\def\epbar{\overline{\ep}}
\def\gammabar{\overline{\gamma}}
\def\omegabar{\overline{\omega}}

\def\rhohat{\hat{\rho}}
\def\tauhat{\hat{\tau}}

\def\nutilde{\tilde{\nu}}

\def\del{{\partial}}
\def\delbar{\overline{\partial}}
\def\deldelbar{\del\delbar}

\def\bs{\backslash}

\def\blank{\phantom{x}}
\def\bil#1#2{\langle #1,#2 \rangle}
\def\metric{{|\blank |}}

\def\Sp{{\mathrm{Sp}}}
\def\GL{{\mathrm{GL}}}
\def\SL{{\mathrm{SL}}}
\def\sing{{\mathrm{sing}}}

\newcommand\id{\operatorname{id}}
\newcommand\Hom{\operatorname{Hom}}

\newcommand\Ext{\operatorname{Ext}}
\newcommand\Aut{\operatorname{Aut}}
\newcommand\Diff{\operatorname{Diff}}
\newcommand\Jac{\operatorname{Jac}}
\newcommand\Pic{\operatorname{Pic}}
\newcommand\Gr{\operatorname{Gr}}

\renewcommand\Im{\operatorname{Im}}
\renewcommand\Re{\operatorname{Re}}

\newtheorem{theorem}{Theorem}[section]
\newtheorem{lemma}[theorem]{Lemma}
\newtheorem{proposition}[theorem]{Proposition}
\newtheorem{corollary}[theorem]{Corollary}
\newtheorem{bigtheorem}{Theorem}
\newtheorem{bigcor}[bigtheorem]{Corollary}

\theoremstyle{definition}
\newtheorem{definition}[theorem]{Definition}

\theoremstyle{remark}
\newtheorem{remark}[theorem]{Remark}

\begin{document}

\title{On the Arakelov Geometry of Moduli Spaces of Curves}

\author{Richard Hain}
\address{Department of Mathematics\\ Duke University\\
Durham, NC 27708-0320}
\email{hain@math.duke.edu}

\author{David Reed}
\address{Department of Mathematics\\ Duke University\\
Durham, NC 27708-0320}
\email{dreed@math.duke.edu}

\date{\today}

\thanks{The first author was supported in part by grants from the
National Science Foundation.}


\keywords{moduli space of curves, heights, Arakelov theory, Faltings'
delta function, Torelli group, mapping class group}

\maketitle

\section{Introduction}

In this paper we consider some problems in the Arakelov geometry of
$\M_g$, the moduli space of smooth projective curves of genus $g$ over
$\C$. Specifically, we are interested in naturally metrized line bundles
over $\M_g$ and their extensions to $\Mbar_g$, the Deligne-Mumford
compactification of $\M_g$. These line bundles typically occur when
computing the archimedean height of a curve. Here they arise in
computations of the archimedean height of the algebraic cycle $C-C^-$
in the jacobian of a genus $g$ curve $C$ (cf.\ \cite{hain:heights}).

As is well known, $\Pic\M_g$ is of rank $1$ when $g\ge 3$. It is
generated by the
determinant bundle $\L:=\det \pi_*\Omega_{\cC/\M_g}^1$, where $\cC$
denotes the universal curve over $\M_g$.  It has a natural metric
which is induced from the standard one on $\pi_*\Omega_{\cC/\M_g}^1$:
$$
\|\omega\|^2 =\frac{i}{2}\int_C\omega \wedge \omegabar,\qquad
\omega \in H^0(C,\Omega^1).
$$
For any metrized line bundle $(\N, \metric_\N)$ over $\M_g$, there is some
integer $N$ such that $\mathcal N \cong \L^{\otimes N}$.  This isomorphism
is unique up to a constant function as there are no non-constant invertible
functions on $\M_g$.  Hence any two metrized line bundles over $\M_g$ can
be ``compared"; there is a real valued function $f: \M_g \to \R$, unique up
to a constant, so that
$$
\metric_\N = e^f \metric_{\L^{\otimes N}}.
$$
From this it follows that, when $g\ge 3$, we have a group isomorphism
$$
\left\{\text{\parbox[c]{1.8in}{isomorphism classes of metrized
line bundles over $\M_g$}}\right\}
\cong
\left\{(N,df) : N\in \Z \text{ and } f : \M_g \to \R\right\}
$$
under which $(\N,\metric_\N)$ corresponds to $(N, df)$. Thus the study of
interesting metrized line bundles over $\M_g$ is just the study of interesting
real valued functions on $\M_g$ (mod constants).

Faltings \cite{faltings} used such a comparison to define functions
$\d_g : \M_g \to \R$ for all $g\ge 1$. He did this by constructing a second
metric on $\L$. Since both metrics are on the same line bundle, the function
$\d_g$ is a genuine function on $\M_g$, and not merely defined mod constants.

In this paper we consider a third naturally metrized line bundle
$(\B, \metric_\B)$ over $\M_g$ for all $g \ge 3$. As a line bundle, $\B$
is isomorphic to the $(8g+4)$th power of $\L$. The line bundle $\B$ is
the {\it biextension line bundle} associated to the algebraic 1-cycle
$C - C^-$ in the jacobian $\Jac C$
of a genus $g$ curve $C$. We review the construction of such biextension
bundles in Section~\ref{biext}.  Denote the logarithm of the ratio of the
metrics on $\B$ and $\L^{\otimes(8g+4)}$ by $\b_g : \M_g \to \R$.  It is
well defined mod constants.

For two functions $f$ and $g$ defined on a punctured disk about $0$, we
write $f_1 \sim f_2$ as $t\to 0$ to mean $f_1(t)-f_2(t)$ remains bounded as
$t\to 0$.

\begin{bigtheorem}
\label{asym}
Suppose that $X \rightarrow \disk$ is a proper family of stable curves of genus
$g \ge 3$ over the unit disk. Suppose that $X$ is smooth and that each
$X_t$ is smooth  when $t\neq 0$.
\begin{enumerate}
\item If $X_0$ is irreducible and has only one node, then
$$
\beta_g(X_t) \sim -g \log|t|- (4g+2)\log\log(1/|t|)\text{ as }t\to 0.
$$
\item If $X_0$ is reducible with one node and its components have genera
$h$ and $g-h$, then
$$
\beta_g(X_t) \sim -4h(g-h)\log|t| \text{ as }t\to 0.
$$
\end{enumerate}
\end{bigtheorem}

These asymptotics may be compared with those of Faltings'
$\delta$-functions which were established by Jorgenson \cite{jorgenson}
and Wentworth \cite{wentworth}.  In the notation of Theorem~\ref{asym},
their result is:

\begin{enumerate}
\item If $X_0$ irreducible with a single node,
$$
3g\delta_g(X_t) \sim -(4g-1)\log|t| - 18 \log\log(1/|t|)\text{ as }t\to 0.
$$
\item If $X_0$ reducible with a single node and components of genera $h$
and $g-h$
$$
3g\delta_g(X_t) \sim -12h(g-h)\log|t|\text{ as }t\to 0.
$$
\end{enumerate}

The incommensurablity of the asymptotics implies:

\begin{bigcor}
If $g\ge 3$, then the exact $1$-forms
$d\d_g$ and $d\b_g$ on $\M_g$ are linearly independent over $\R$.
\end{bigcor}

An explicit formula for $\b_3$ is given in \cite{hain-reed:formula}. It
involves the Siegel modular form $\chi_{18}$ and an integral. Although
explicit, the integral does not seem to be computable at this time.

Motivated by the formula for $\b_3$, and the asymptotics of $\b_g$ when
$g\ge 3$ it seems natural to define
$$
\b_1 = -\log \|\D\| \text{ and }\b_2 = -2\log \|\chi_{10}\|
$$
where
$$
\D(q) = (2\pi)^{12}q\prod_{n\ge 0}(1-q^n)^{24}
$$
is the cusp form of weight 12 in genus 1,
$$
\chi_{10}(\Omega) = \prod_{\alpha \text{ even}}\vartheta_\alpha(0;\Omega)^2
$$
is the Siegel modular form of weight 10 on $\h_2$, the Siegel upper half
plane of rank 2, and $\|\blank\|$ denotes the norm
$$
\|F(\Omega)\| = |F(\Omega)|(\det \Im \Omega)^{k/2},\qquad \Omega \in \h_g
$$
of a modular form of weight $k$.
Both $\b_1$ and $\b_2$ satisfy the asymptotic formulas in Theorem~\ref{asym}.
It would be desirable to have a uniform construction of the $\b_g$ for all
$g\ge 1$ that produces the $\b_1$ and $\b_2$ defined above.

Denote the Deligne-Mumford compactification of $\M_g$ by $\Mbar_g$ and the
components of the boundary divisor $\Mbar_g - \M_g$ by $\D_h$, where $0 \le h
\le [g/2]$. The (orbifold) singular locus $\D_0^\sing$ of $\D_0$ is the closure
of the locus of stable genus $g$ curves $C$ with two nodes, neither of which
will disconnect $C$ when removed.

Theorem~\ref{asym} is a direct consequence of the following Theorem and a
result of Dale Lear \cite{lear}, as explained in Section~\ref{overview}.

\begin{bigtheorem}
\label{main}
The biextension line bundle $\B$ over $\M_g$ extends naturally to a line bundle
$\Bbar$ on $\Mbar_g$. The metric extends to a smooth metric on $\Bbar$ over
$\Mbar_g - \D_0$ and extends continuously to the restriction of $\Bbar$ to any
holomorphic arc meeting $\D_0$ transversally. The first Chern class of the
extended bundle is
$$
(8g+4)\lambda - g \d_0 - 4 \sum_{h=1}^{[g/2]} h(g-h)\d_h.
$$
\end{bigtheorem}

At present we are not sure whether the metric extends continuously to all of
$\Mbar_g$. This will probably have to wait until our understanding of limits
of mixed Hodge structures improves. Work in this direction is being done by
Gregory Pearlstein \cite{pearlstein} and we hope that this question will be 
resolved in his work.

Also, it is interesting to note that the extended line bundle $\Bbar$ is the
line bundle shown by Moriwaki \cite{moriwaki} to have non-negative degree on
every complete curve in $\Mbar_g$ not contained in the boundary divisor $\D$.
This is no accident. In a subsequent paper, the first author will show that the
curvature of $\Bbar$ on $\Mbar_g - \D_0$ is a non-negative 2-form which is
locally $L_1$ on every holomorphic arc intersecting $\D_0$, thus giving an
analytic proof of Moriwaki's inequality.

\section{Preliminaries}

All varieties in this paper will be over $\C$. The moduli space of $n$-pointed
($n\ge 0$) smooth projective curves will be denoted by $\M_g^n$ and its
Deligne-Mumford compactification by $\Mbar_g^n$. We will be primarily
interested in the cases when $n=0$ and $1$.  When $n=0$, we will write $\M_g$
and $\Mbar_g$ instead of $\M_g^0$ and $\Mbar_g^0$.

The moduli space of principally polarized abelian varieties of dimension $g$
will be denoted by $\A_g$. The level $l$ covers of $\M_g$ and $\A_g$ will be
denoted by $\M_g[l]$ and $\A_g[l]$.

All moduli spaces will be regarded as orbifolds, or more accurately, as stacks
in the sense of Mumford \cite{mumford}. By the fundamental group we shall mean
the orbifold fundamental group. In particular, when $g\ge 2$, the fundamental
group of $\M_g^n$ is isomorphic to $\G_g^n$, the mapping class group of an
$n$-pointed,  genus $g$ surface. The fundamental group of $\A_g$ is isomorphic
to $\Sp_g(\Z)$. 

The difference $\Mbar_g^n - \M_g$ is a divisor $\D$ which is smooth
with normal crossings (provided we regard $\Mbar_g^n$ as an orbifold).
When $n=0$
$$
\D = \bigcup_{h=0}^{[g/2]} \D_h
$$
The generic point of each each $\D_h$ corresponds to a stable genus $g$ curve
with 1 node. The generic point of $\D_0$ is irreducible and has normalization
of genus $g-1$; the generic point of $\D_h$ when $h>0$ corresponds to a
reducible curve with components of genera $h$ and $g-h$.

We shall denote by $\Mtilde_g$ the open subset $\Mbar_g - \D_0$. The period map
$\M_g \to \A_g$ extends to $\Mtilde_g$ and the extended period map $\Mtilde_g
\to \A_g$ is proper and has closed image. We can also define $\Mtilde_g^n$ to
be the inverse image of $\Mtilde_g$ under the  natural map $\Mbar_g^n \to
\Mbar_g$.

The first Chern class of the determinant of the Hodge bundle will be denoted by
$\lambda$. Since $\Mbar_g$ is a compact orbifold, the divisor $\D_h$ determines
a class $\d_g \in H^2(\Mbar_g,\Q)$.  Recall from \cite{arb-corn} that
$\Pic\Mbar_g$ is free with basis $\lambda, \d_0, \d_1, \dots , \d_{[g/2]}$.
From this it follows that $\Pic \Mtilde_g$ is also free with basis $\lambda, 
\d_1, \dots , \d_{[g/2]}$

The Torelli group $T_g^n$ is the kernel of the natural map $\G_g^n \to
\Sp_g(\Z)$ induced by the period map $\M_g^n \to \A_g$.

Denote the moduli space of smooth hyperelliptic curves of genus $g$ by $\H_g$.
We view it as an orbifold. The following lemma is surely well known, but we
state and prove it for lack of a convenient reference.

\begin{lemma}\label{constants}
For all $g\ge 1$ and all $n\ge 0$, the only invertible functions on $\M_g^n$,
$\H_g$ and $\A_g$ are the constants. That is,
$$
H^0(\M_g^n,\O^\times) = H^0(\H_g,\O^\times) = H^0(\A_g,\O^\times) = \C^\times.
$$
\end{lemma}

\begin{proof}
Suppose that $X$ is a smooth variety. If $f$ is a non-constant invertible
function on $X$, then $d\log f/2\pi i$ is a non-trivial element of $H^1(X,\Q)$.
Each of the moduli spaces considered in this proposition has trivial first
homology with rational coefficients. This is well known for $\M_g$ and $\A_g$.
The corresponding statement for $\H_g$ follows as $\H_g$ is the moduli space of
$(2g-2)$-tuples of distinct points in $\P^1$ modulo projective equivalence:
$$
\H_g = \big((\C-\{0,1\})^{2g-1} - \Delta\big)/S_{2g+2},
$$
where $\Delta$ denotes the fat diagonal. Since the first homology of
$(\C-\{0,1\})^{2g-1}-\Delta$ does not contain a copy of the trivial
representation of $S_{2g+2}$, it follows that $H_1(\H_g,\Q)$ is
trivial.\footnote{The first homology is the irreducible $S_{2g+2}$-module
corresponding to the partition $[2g,2]$ of $2g+2$. A proof can be found in
\cite[(10.6)]{hain-macp}. There our space is denoted $Y^2_{2g-1}$, which is
$Y^{2g}_1$ by duality --- see \cite[(5.7)]{hain-macp}.} It follows that none of
$\M_g$, $\A_g$ or $\H_g$ has a non-constant invertible function when $g\ge 1$.
\end{proof}

Finally, some notation. If $A$ is a $\Z$ (mixed) Hodge structure, 
we shall use $A_\Z$, $A_\R$, $A_\C$ to denote the integral lattice,
the underlying real vector space, and the underlying complex vector space.
If $A$ has negative weight (typically $-1$), we shall set
$$
J(A) = A_\C/(F^0A + A_\Z).
$$
When $A_\Z$ is torsion free, this is naturally isomorphic to $\Ext^1_\H(\Z,A)$,
where $\H$ denotes the category of mixed Hodge structures.

\section{Some Linear Algebra}
\label{lin_alg}

In this section we dispense with some linear algebra. Suppose $g \ge 2$.
Fix a lattice $H$ of
rank $2g$ with a unimodular symplectic form $(\blank \cdot \blank )$. This 
gives a natural identification of $H$ with its dual. For example, we will
take $H$ to be $H_1(C,\Z)$ where $C$ is a compact Riemann surface with the
intersection paring, or more generally, $H_1(A,\Z)$ where $A$ is a principally
polarized abelian variety. The dual of the intersection pairing may be viewed
as an element $\theta$ of $\Lambda^2 H$. Set 
$$
V = \Lambda^3 H/ (\theta \wedge H).
$$
Denote the group of automorphisms of $H$ that preserve the intersection
pairing by $\Sp(H)$. By choosing a symplectic basis of $H$, we can identify
it with $\Sp_g(\Z)$. Note that $V$ is an $\Sp(H)$-module. 

Denote the intersection number of $x,y\in H$ by $x\cdot y$. This can
be used to define an $\Sp(H)$-equivariant map
$$
c : \Lambda^3 H \to H.
$$
It takes $x \wedge y \wedge z$ to $(x\cdot y)z + (y\cdot z)x + (z\cdot x)y$.
There is also the $\Sp(H)$-invariant map $H \to \Lambda^3 H$ that
takes $x$ to $x\wedge \theta$.

\begin{proposition}
\label{comp}
For all $x\in H$, $c(x\wedge \theta) = (g-1)x$. Consequently, if $g = 2$,
then $c : \Lambda^3 H \to H$ is an isomorphism and $V=0$. \qed
\end{proposition}

This shows that
$$
\Lambda^3 H_\Q \cong V_\Q \oplus H_\Q
$$
as $\Sp(H)$ modules. It is well known from the representation theory of the
symplectic group that both $V_\Q$ and $H_\Q$ are irreducible as $\Sp(H)$-%
modules and that each admits an invariant symplectic form, unique up to
a constant multiple. Next we give a formula for the primitive symplectic
form $q : \Lambda^2 V \to \Z$.

Denote the image of $x \wedge y \wedge z$ in $V$ by 
$
\overline{x \wedge y \wedge z}.
$
It follows from the above that the function
$$
j : V \to \Lambda^3 H
$$
defined by
$$
\overline{x\wedge y\wedge z} \to
(g-1)x\wedge y \wedge z - \theta\wedge c(x\wedge y\wedge z)
$$
is integral and $\Sp_g(\Z)$ invariant.
Denote the symplectic form on $\Lambda^3 H$ induced from that
on $H$ by $\bil \blank \blank$. That is,
$$
\bil{x_1\wedge x_2 \wedge x_3}{ y_1\wedge y_2\wedge y_3} = \det(x_i\cdot y_j)
$$

The symplectic form $q$ on $V$ is defined by
$$
q(u,v) = \bil{j(u)}{j(v)} /(g-1).
$$
An easy computation shows that
$$
\bil{x\wedge \theta} {y\wedge \theta} = (g-1)(x\cdot y)
$$
Since $\theta\wedge H$ and $j(V)$ are orthogonal, it follows that $q$ takes
values in $\Z$.

\section{The Fundamental Normal Function}

In this section we recall the construction of the fundamental normal
function over $\M_g$. More details can be found in \cite{hain:normal}
and \cite[\S 7]{hain-looijenga}.

Suppose that $A$ is a principally polarized abelian variety of dimension
$g$. Let
$$
V_{A,\Z} = \Lambda^3 H_1(A,\Z) / H_1(A,\Z) = H_3(A,\Z) / H_1(A,\Z)
$$
Here $H_1(A,\Z)$ is imbedded in $\Lambda^3 H_1(A,\Z)$ via the map
that takes $x \in H_1(A,\Z)$ to $x \wedge \theta$, where
$$
\theta\in \Lambda^2 H_1(A,\Z)
$$
is the dual of the principal polarization.
Note that when $g \le 2$, $V_{A,\Z}$ is trivial --- cf.\ Proposition~\ref{comp}.

There is unique Hodge structure of weight $-1$ on $V_A$ such that
the quotient map
$$
H_3(A,\Z(-1)) \to V_A
$$
is a morphism of Hodge structures. From this we can construct the
compact complex torus
$$
J(V_A) := V_{A,\C} / \left( F^0 V_A + V_{A,\Z}\right).
$$
This can be thought of as the Griffiths intermediate jacobian associated
to the primitive part of the third homology of $A$.

When $A$ is the jacobian $\Jac C$ of a smooth projective curve $C$
of genus $g$, we shall abuse notation and write $V_C$ instead of
$V_{\Jac C}$.

Harris's {\it Harmonic volume} \cite{harris, pulte} gives a canonical
point $\nu(C)$ of $J(V_C)$ whose construction we now recall.
Choose a point $x \in C$. This gives an imbedding
$$
C \to \Jac C
$$
by taking $y\in C$ to $[y] - [x]$. Denote the image of this map by 
$C_x$ and the image of $C_x$ under the involution $D \mapsto -D$ of $\Jac C$
by $C_x^-$. The cycle $C_x - C_x^-$ is homologous to 0, and therefore
gives a point $\nu_x(C)$ in the Griffiths intermediate jacobian
$$
J_1(\Jac C) :=
H_3(\Jac C,\C)/\left(F^{-1} H_3(\Jac C) + H_3(\Jac C,\Z)\right).
$$
This compact complex torus maps to $J(V_C)$ as $H_3(\Jac C)$ is  canonically
isomorphic to $\Lambda^3 H_1(C,\Z)$. The kernel of the quotient map $J(\Jac C)
\to J(V_C)$ is naturally isomorphic to $\Jac C$. One has \cite{pulte} that
$$
\nu_x(C) - \nu_y(C) = - 2([x]-[y]) \in \Jac C.
$$
It follows that the image of $\nu_x(C)$ in $J(V_C)$ is independent
of $x$. The common image of the $\nu_x(C)$ in $J(V_C)$ will be denoted by
$\nu(C)$.

The abelian groups $V_{A,\Z}$ form a local system, or more accurately, a 
polarized variation of Hodge structure of weight $-1$ over $\A_g$. We shall
denote it by $\V$. Taking the intermediate jacobian of each fiber gives a
complex analytic bundle
$$
\J(\V) \to \A_g
$$
whose fiber over the moduli point $[A]$ of $A$ is $J(V_A)$.

\begin{remark}\label{low_genus}
This bundle of tori is trivial (i.e., the fiber is 0) when $g \le 3$.
\end{remark}

Harmonic volume gives a lift of the period map $\M_g \to \A_g$:
$$
\xymatrix{
 & {\J(\V)} \ar[d] \cr
{\M_g} \ar[ur]^{\nu} \ar[r] & {\A_g}
}
$$

We shall denote the variety $\Mbar_g - \D_0$ by $\Mtilde_g$. This variety
is the moduli space of stable genus $g$ curves whose dual graph is a tree,
or equivalently, whose $\Pic^0$ is an abelian variety. The period map
extends to a proper map $\Mtilde_g \to \A_g$.

\begin{proposition}
The normal function $\nu$ extends to  a normal function
$\nutilde : \Mtilde_g \to \J(\V)$ that lifts the period mapping
$\Mtilde_g \to \A_g$.
\end{proposition}

\begin{proof}
This follows immediately from \cite[(7.1)]{hain:normal} as the variation
of Hodge structure $\V$ extends to a variation over $\Mtilde_g$.
\end{proof}

\section{Chern Classes}

Suppose that $T \to B$ is a smooth family of compact tori over a connected
manifold $B$. Fix a base point $b_o$ of $B$ and denote the fiber of it by
$F_o$. We shall also assume that the family has a distinguished section,
which we shall call the zero section.

\begin{lemma}
\label{gen_lifts}
For such a family of compact tori, there is a natural mapping
$$
s : H^0(B,H^k(F_o,\R)) \to H^k(T,\R)
$$
whose composition with the projection
$$
H^k(T,\R) \to H^0(B,H^k(F_o,\R))
$$
is the identity. Moreover, for each $u\in H^0(B,H^k(F_o))$, the extended class
$\utilde := s(u)$ has a natural representative $\wtilde_u$ whose restriction to
the zero-section is trivial and to each fiber is translation invariant.
\end{lemma}

\begin{proof}
Denote the fiber of $E \to B$ over $b$ by $F_b$. We have the local systems
$\V_\Z$ and $\V_\R$ whose fibers over  $b \in B$ are $H_1(F_b,\Z)$ and
$H_1(F_b,\R)$. The quotient $\V_\R/\V_\Z$ is a family of tori over $B$
which is naturally isomorphic to the family $T\to B$ via an isomorphism that
preserves the zero sections.

Denote the base point of $B$ by $b_o$. The first homology of $F_o$, the
fiber over $b_o$ is a $\pi_1(B,b_o)$-module. Note that the group
$$
\pi_1(B,b_o) \ltimes H_1(F_o,\Z)
$$
acts on $\widetilde{B} \times H_1(F_o,\R)$ via the action
$$
(\gamma, z) : (b,a) \mapsto (\gamma\cdot b, \gamma\cdot a + z).
$$
Here $\widetilde{B}$ denotes a universal covering of $B$. The quotient
$$
\left(\pi_1(B,b_o) \ltimes H_1(F_o,\Z)\right)
\bs
\left(\widetilde{B} \times H_1(F_o,\R)\right)
$$
is naturally isomorphic to $\V_\R/\V_\Z$ as a bundle over $B$.

Each $\pi_1(B,b_o)$-invariant cohomology class $u$ in $H^k(F_o,\R)$ is
represented by a translation invariant $k$-form $\w_u$ on $H_1(F_o,\R)$.
The pull back of $\w_u$ along the projection
$p:\widetilde{B} \times H_1(F_o,\R) \to H_1(F_o,\R)$ is a $k$-form
invariant under the action of $\pi_1(B,b_o) \ltimes H_1(F_o,\Z)$.
It therefore descends to a closed form $\wtilde_u$ on $\V_\R/\V_\Z$ whose
restriction to the fiber over $b_o$ is $\w_u$.

The mapping $s$ is defined by taking $s(u)$ to be the de~Rham class of
$\wtilde_u$.
\end{proof}

\begin{corollary}
With the same hypotheses as above, for all $r \ge 2$, the differentials
$$
d_r : E_r^{0,k} \to E_r^{r,k-r+1}
$$
vanish in the Leray-Serre spectral sequence
$$
E_2^{s,t} = H^s(B,H^t(F_o,\R)) \implies H^{s+t}(T,\R). \qed
$$
\end{corollary}

\begin{corollary}
\label{lift}
If $H^1(B,H^1(F_o,\R))$ vanishes, then each $\pi_1(B,b_o)$-invariant element
$u$ of  $H^2(F_o,\R)$ extends to a unique element $\utilde$ of $H^2(T,\R)$
whose restriction to $F_o$ is $u$ and whose restriction to the zero section is
trivial. If $u$ is a rational cohomology class, then so is $\utilde$.
\end{corollary}

\begin{proof}
The invariant class $u$ lifts to $\utilde = s(u) \in H^2(T,\R)$, which has the
property that its restriction to the zero section vanishes. The vanishing of
$H^1(B,H^1(F_o,\R))$, Lemma~\ref{gen_lifts} and the existence of a section (the
zero section) imply that the sequence
$$
0 \to H^2(B) \to H^2(T) \to H^2(F_o)^{\pi_1(B,b_o)} \to 0
$$
is exact with real, and therefore rational, coefficients. Restricting to the
zero section defines a splitting of the left hand mapping. This gives a direct
sum decomposition of $H^2(T)$ which is defined over $\Q$. The uniqueness and
rationality statements follow.
\end{proof}

\begin{corollary}
\label{unique}
For all $g \ge 1$, there is a canonical element $\phi$ of $H^2(\J(\V),\Q)$
whose restriction to the zero section is trivial and whose restriction to each
fiber is the class of the polarization $q$. It is characterized by these
properties.
\end{corollary}

\begin{proof}
When $g\le 2$, the the result is trivially true as the bundle $\J(\V)$ has the
trivial complex torus as fiber. When $g\ge 2$, we know that $H^1(\Sp_g(\Z),V)$
vanishes \cite{raghunathan}, and so one can apply the previous result.
\end{proof}

We will see in Proposition~\ref{integral} that $2\phi$ is integral. The
following result is due to Morita \cite[(5.8)]{morita:chern}. The precise
statement below and a more direct proof can be found in \cite{hain-reed:chern}.

\begin{theorem}[Morita]
\label{morita}
For all $g\ge 1$, we have $\nu^\ast(2\phi) = (8g+4)\lambda$ in $H^2(\M_g,\Z)$.
\end{theorem}

\section {A Characterization of the Biextension Line Bundle and its Metric}
\label{biext}

Since the construction of the biextension line bundle and its metric
relies on abstract Hodge theory and may not be immediately accessible, we
first give a characterization of it and its metric. The actual
construction is reviewed in the next section.

\begin{proposition}
\label{chern_char}
There is at most one holomorphic line bundle $\Bhat$ over $\J(\V)$ whose
first Chern class is $2\phi$ and whose restriction to the zero
section is trivial.
\end{proposition}

\begin{proof}
Two such bundles differ by a line bundle which is trivial on
the zero section and is topologically trivial on each fiber. The group
of topologically trivial line bundles on $J(V_A)$ is $J(\Vdual_A)$
where
$$
\Vdual = \Hom(V,\Z(1))
$$
(cf.\ \cite[(3.1.6)]{hain:heights}).
Thus, to show that two such bundles are isomorphic, it suffices
to show that the bundle of intermediate jacobians associated to the
variation $\check{\V}$ has no sections. It follows from
\cite[(9.2)]{hain:normal} that all normal function sections
of $\J(\check{\V})$ are torsion. The result therefore follows from
the following algebraic fact.
\end{proof}

\begin{proposition}
For each prime $p$, there are no non-zero $\Sp_g(\Z)$-invariant elements of
$\Hom_\Z(V_\Z,\F_p)$, $\Hom_\Z(\Lambda^3 H_\Z,\F_p)$ or $\big(\Lambda^3
H_\Z\big)\otimes \F_2$.
\end{proposition}

\begin{proof}
Note that
$$
\Hom_\Z(V_\Z,\F_p) \cong \Hom_{\F_p}(V_\Z\otimes \F_p,\F_p).
$$
By the right exactness of tensor product, the natural mapping
$$
\big(\Lambda^3 H_\Z\big)\otimes \F_p \to V_\Z \otimes \F_p
$$
is surjective, from which it follows that $\Hom_\Z(V_\Z,\F_p)$ is
a submodule of
$$
\Hom_{\F_p}(\big(\Lambda^3 H_\Z\big)\otimes \F_p,\F_p).
$$
Thus it suffices to prove that
$$
\big[\big(\Lambda^3 H_\Z\big)\otimes \F_2\big]^{\Sp_g(\Z)} = 0 \text{ and }
\Hom_{\Sp_g(\Z)}\big(\big(\Lambda^3 H_\Z\big)\otimes \F_2,\F_2\big) = 0.
$$

When $p$ is odd, $-I \in \Sp_g(\Z)$ acts as $-1$ on $\Hom_\Z(\Lambda^3
H_\Z,\F_p)$ and $\big(\Lambda^3 H_\Z\big)\otimes \F_2$, which implies there are
no invariants except possibly when $p=2$.

Choose a symplectic basis $a_1,\dots,a_g,b_1,\dots,b_g$ of $H_\Z$.
Let $H_j$ be the span of $a_j$ and $b_j$. Then
$$
H_\Z = H_1 \oplus \dots \oplus H_g.
$$
The intersection of $\Sp_g(\Z)$ with $\GL(H_j)$ is isomorphic to $\SL_2(\Z)$
and this gives an imbedding of $\SL_2(\Z) \times \dots \times \SL_2(\Z)$
($g$ factors) into $\Sp_g(\Z)$.

Denote the mod $2$ reduction of $H_\Z$ by $W$ and of $H_j$ by $W_j$. Set
$$
U = \big(\Lambda ^3 H_\Z\big)\otimes \F_2
$$
Since $\SL_2(\F_2)$ is isomorphic to $S_3$, the symmetric group on 3 letters,
we have a diagonal imbedding of $G = S_3 \times \dots \times S_3$ ($g$ factors)
into the image of $\Sp_g(\Z)$ in $\GL(U)$.

Since $(\F_2)^2$ is a {\it simple} $S_3$-module, it follows that each
$W_j$ is a simple $G$-module. This implies that if $i<j<k$, then
$W_i \otimes W_j \otimes W_k$ is also a {\it simple} $G$-module. Since
$$
U \cong
\bigoplus_{j\neq k} W_j \otimes \Lambda^2 W_k \oplus 
\bigoplus_{i<j<k} W_i \otimes W_j \otimes W_k
\cong (W_1 \oplus \dots \oplus W_g)^{g-1} \oplus 
\bigoplus_{i<j<k} W_i \otimes W_j \otimes W_k,
$$
it is also a direct sum of simple $G$-modules, none of which is trivial. It
follows that $U^G = 0$ and $\Hom_G(U,\F_2) = 0$. The result follows.
\end{proof}

Since the group of $p$-torsion sections of $\J(\bVdual)$ defined over $\M_g$
is isomorphic to $\Hom_\Z(V_\Z,\F_p)^{\Sp_g(\Z)}$, we have:

\begin{corollary}
There are no non-trivial torsion sections of $\J(\bVdual)$ defined over $\M_g$.
\end{corollary}

Since the only invertible functions on $\A_g$ are constants
(Lemma~\ref{constants}), any two  trivializations of a trivial line bundle over
$\A_g$ differ by a constant. In particular, if the line bundle $\Bhat$ over
$\J(\V)$ exists, the trivialization of its restriction to the zero section is
unique up to a constant. In particular, it makes sense to say that a metric on
the restriction of $\Bhat$ to the zero section is constant.

In the next section, we will construct the line bundle $\Bhat$ over $\J(\V)$
whose first Chern class is $2\phi \in H^2(\J(\V))$. We also construct a metric
$\metric_\Bhat$ on it. The metric has the property that the  restriction of its
curvature to each fiber is translation invariant and to the zero section
vanishes. Proposition~\ref{chern_char} and the following result show that these
properties characterize $\Bhat$ and its metric.

\begin{proposition}
If such a line bundle $\Bhat$ over $\J(\V)$ exists, then any two metrics
on it that satisfy:
\begin{enumerate}
\item the metric on the restriction of $\Bhat$ to the zero section
is constant;
\item the curvature of the restriction of the metric to each fiber
is translation invariant;
\end{enumerate}
are constant multiples of each other.
\end{proposition}

\begin{proof}
The second condition determines the metric on each fiber up to a constant
multiple. The first condition then determines all the fiber metrics up to
a common, constant multiple.
\end{proof}

\begin{remark}
\label{low_genus2}
Note that the since the bundle $\J(\V)$ is trivial when $g = 2$, the line
bundle $\Bhat$ is also trivial in that case and has the trivial metric (cf.\
Remark~\ref{low_genus}).
\end{remark}

\begin{definition}
The biextension line bundle $\B \to \Mtilde_g$ is the pullback of the
biextension line bundle $\Bhat \to \J(\V)$ along the fundamental normal
function $\nutilde : \Mtilde_g \to \J(\V)$. Its metric $\metric_\B$ is the
pullback of $\metric_\Bhat$ along $\nutilde$.
\end{definition}

We will use the same notation to denote the restriction of the biextension
bundle to $\M_g$. It follows from Morita's Theorem (Thm.~\ref{morita}) that the
first Chern class of the biextension bundle $\B \to \M_g$ is $(8g+4)\lambda$.

Note that it makes sense to talk about the constant metric on the trivial
bundle over the hyperelliptic locus $\H_g$ as the only invertible functions on
$\H_g$ are constants (cf.\ Lemma~\ref{constants}).

\begin{proposition}
\label{hyperelliptic}
The restriction of the biextension line bundle to the hyperelliptic locus
is the trivial bundle with a constant metric. Consequently, the restriction
of the curvature of $(\B,\metric_\B)$ to $\H_g$ vanishes. \qed
\end{proposition}

\begin{proof}
Since the cycle $C_x - C_x^-$ vanishes identically when $C$ is hyperelliptic
and $x$ is a Weierstrass point, it follows that the normal function $\nu$
vanishes identically on the hyperelliptic locus $\H_g$. The result follows
as the biextension bundle and its metric are trivial on the zero-section.
\end{proof}

\section{Construction of the Biextension Bundle}

In this section we review the construction of the biextension line bundle
$\Bhat \to \J(\V)$ and its metric $\metric_\Bhat$ (cf. \cite{hain:heights},
and \cite[\S8]{hain:comp}).

We begin in the more general situation of a variation of Hodge structure $\U$
of weight $-1$ over a complex manifold $X$. Denote the dual variation
$\Hom(\U,\Z(1))$ by $\bUdual$; it is also of weight $-1$. The fibers of $\U$
and $\bUdual$ over $x\in X$ will be denoted by $U_x$ and $\Udual_x$,
respectively.

One has the corresponding bundles of intermediate jacobians $\J(\U)$ and
$\J(\bUdual)$ over $X$. Their fibers $J(U_x)$ and $J(\Udual_x)$ over $x \in X$
are dual complex tori \cite[(3.1.6)]{hain:heights} and can be naturally
identified with $\Ext_\H(\Z,U_x)$ and $\Ext_\H(U_x,\Z(1))$, respectively. For
this reason, we shall denote the family of intermediate jacobians $\J(\bUdual)$
by $\Jdual(\U)$.

One can form the bundle
$$
\J(\U)\times_X \Jdual(\U) \to X.
$$
There is a natural principal $\C^\ast$-bundle $\B(\U)^\ast$ over this. Its
fiber $\B(\U)^\ast_x$ over $x\in X$ is the set of mixed Hodge structures $B$
with 
$$
\Gr^W_k B = 0 \text{ if $k \neq 0,-1,-2$}
$$
together with isomorphisms
$$
\Gr^W_{0} B \cong \Z,\,
\Gr^W_{-1} B \cong U_x \text{ and } \Gr^W_{-2} B \cong \Z(1).
$$
The projection
$$
\B(\U)^\ast \to \J(\U)\times_X \Jdual(\U)
$$
takes $B$ to $(B/\Z(1),W_{-1}B)$. There is a natural action of $\C^\ast$ on the
fibers which is described in \cite[(3.2.2)]{hain:heights}. This $\C^\ast$ can
be naturally identified with $\Ext_\H^1(\Z,\Z(1))$.

Denote the corresponding line bundle by $\B(\U)$. It has a natural metric which
is described in \cite[(3.2.9)]{hain:heights} and can be described briefly as
follows: The moduli space of {\it real} mixed Hodge structures with weight
graded quotients $\R, U_{x,\R}$ and $\R(1)$ is naturally isomorphic to $\R$.
Taking $B$ to the canonical real period of $B\otimes \R$ gives a function
$\mu : \B(\U)^\ast \to \R$. There is an explicit formula
\cite[(3.2.11)]{hain:heights} for $\mu(B)$ in terms of the $\Z$-periods of $B$.
The metric on $\B(\U)^\ast$ is defined by
$$
|B|_\B = \exp \mu(B).
$$

Because of the natural isomorphism
$$
H_1(J(U_x)\times J(\Udual_x)) \cong U_{x,\Z} \oplus \Udual_{x,\Z},
$$
we can canonically identify $H^2(J(U_x)\times J(\Udual_x),\C)$ with the set of
skew symmetric bilinear mappings
$$
\big(U_x \oplus \Udual_x\big)^{\otimes 2} \to \C.
$$
In particular, we have the two dimensional cohomology class $Q$ defined by
$$
Q : \big((u,\psi),(v,\phi)\big) \mapsto \frac{1}{2\pi i} (\psi(v) - \phi(u))
$$
which is integral and invariant under monodromy. It therefore determines a 
closed 2-form $\w_Q$ on $\J(\U)\times_X \Jdual(\U)$ as in
Proposition~\ref{gen_lifts}.

\begin{proposition}
\label{chern_biext}
The first Chern class of the line bundle
$$
\B(\U) \to \J(\U)\times_X \Jdual(\U)
$$
is the class $\w_Q$. It has the property that its restriction to each fiber is
the class corresponding to $Q$ and to the zero section is trivial; these
characterize it when $H^1(X,\U)$ vanishes.
\end{proposition}

\begin{proof}
First note that the restriction of $\B(\U)^\ast$ to the zero section is trivial
as it has the section
$$
x \mapsto \Z \oplus V_x \oplus \Z(1).
$$
It follows from the formula \cite[(3.2.11)]{hain:heights} that this section
has constant length 1.

For the rest of the computations we will use the notation and terminology of 
\cite[\S3]{hain:heights}.  We will consider the fiber of $\B(\U)^\ast$ over 
$x \in X$. Set $U = U_x$ and $\Udual = \Udual_x$. Set
$$
G_{U\, \Z} =
\begin{pmatrix}
1 & U_\Z & \Z(1) \cr 0 & 1 & \Udual_\Z \cr 0 & 0 & 1 \cr
\end{pmatrix}
,\,
G_U =
\begin{pmatrix}
1 & U_\C & \C \cr 0 & 1 & \Udual_\C \cr 0 & 0 & 1 \cr
\end{pmatrix}
,\,
F^0 G_U =
\begin{pmatrix}
1 & F^0 U & 0 \cr 0 & 1 & F^0\Udual \cr 0 & 0 & 1 \cr
\end{pmatrix}
$$
Then the fiber of $\B(\U)^\ast$ over $x\in X$ is
$$
\B(\U)_x^\ast := G_{U\, \Z} \bs G_U /F^0G_U.
$$
The natural projection of this to $J(U) \times J(\Udual)$ is induced by
the obvious group homomorphism $G_U \to U \times \Udual$. We can give
local holomorphic framing $s$ of $\B_x \to J(U) \times J(\Udual)$ by giving
a global holomorphic section of 
$$
G_U /F^0G_U \to U/F^0 U \times \Vdual/F^0\Udual.
$$
We do this by identifying $U/F^0 U$ with $\Fbar^0 U$ and $\Udual/F^0\Udual$
with $\Fbar^0\Vdual$ and then defining
$$
s(g,\gamma) =
\begin{pmatrix}
1 & g & 0 \cr 0 & 1 & \gamma \cr 0 & 0 & 1 \cr
\end{pmatrix}
F^0G_U
$$
where $g \in \Fbar^0 U$ and $\gamma \in \Fbar^0\Udual$.
Using the formula \cite[(3.2.11)]{hain:heights} with 
$$
r = g + \gbar,\, f=-\gbar,\, \rho = \gamma + \gammabar
\text{ and } \psi = -\gammabar
$$
one sees that $\mu(s(g,\gamma)) = - \Re \gamma(\gbar)$.\footnote{Here and
elsewhere in this proof, the conjugate $\hbar \in \Udual$ of $f \in \Udual$ is
the function $\fbar : U \to \C$ defined by $\fbar(v) =
\overline{f(\vbar)}\text{ for all } v \in V$.} It follows that 
$$
\log |s|_\B^2 = -2 \Re \gamma(\gbar).
$$
Choose a basis $e_j$ of $F^0 U$ and a basis $\ep_j$ of $F^0\Udual$ such that
$\ep_j(\ebar_k) = 2\pi i\d_{jk}$ all $j$ and $k$. Since $F^0\Udual$ is the
annihilator of $F^0 U$, $\ep_j(e_k) = 0$ for all $j$ and $k$. Note that
$\epbar_j(e_k) = -2\pi i \d_{jk}$.

Now, if we write $g = \sum_j z_j \ebar_j$ and $\gamma = \sum_j u_j \epbar_j$,
then
$$
\log |s|_\B^2 = - (\gamma + \gammabar)(g + \gbar) =
 2\pi i\, \sum_{j=1}^n \left(z_j \ubar_j - u_j \zbar_j\right).
$$
The Chern class $c_1(\B(U)_x)$ is therefore represented by
$$
\frac{1}{2\pi i}\, \deldelbar \log |s|_\B^2 =
\sum_{j=1}^n \left(dz_j \wedge d\ubar_j - du_j\wedge d\zbar_j\right).
$$
It remains to prove that $c_1(\B(\U))$ is given by $Q$. To see this, we
use the identifications
$$
U_\R \stackrel{\simeq}{\longrightarrow} U_\C/F^0 U
\stackrel{\simeq}{\longleftarrow} \Fbar^0 U
\text{ and }
\Udual_\R \stackrel{\simeq}{\longrightarrow} \Udual_\C/F^0 \Udual
\stackrel{\simeq}{\longleftarrow} \Fbar^0 \Udual.
$$
Set
$$
f_j = e_j + \ebar_j,\, Jf_j = i(e_j - \ebar_j),\, \phi_j = \epbar_j - \ep_j,\,
J\phi_j = i(\epbar_j + \ep_j).
$$
Then the $f_j$ and $Jf_j$ form a basis of $U_\R$, and the $\phi_j$ and
$J\phi_j$ form a basis of $\Udual_\R$. Moreover $f_j$ corresponds to $\ebar_j$,
$\phi_j$ corresponds to $\epbar_j$ and $J$ to multiplication by $i$ under the
isomorphisms above.

Using the standard identification between $U_\R$ and the translation invariant
vector fields on $U_\R/U_\Z$, we have
$$
\langle dz_j,f_k\rangle = \langle d\zbar_j,f_k\rangle = \d_{jk}
\text{ and }
\langle dz_j,Jf_k\rangle = - \langle d\zbar_j,Jf_k\rangle = i\d_{jk}.
$$
Similarly,
$$
\langle du_j,\phi_k\rangle = \langle d\ubar_j,\phi_k\rangle = \d_{jk}
\text{ and }
\langle du_j,J\phi_k\rangle = - \langle d\ubar_j,J\phi_k\rangle = i\d_{jk}.
$$
The formula for the Chern class now follows by direct computation. The final
statement follows directly from Corollary~\ref{lift}.
\end{proof}

We now return to the case where $X = \A_g$ and $\U$ is the variation $\V$. For
each $[A]\in \A_g$, the polarization $q$ gives a MHS morphism $V_A \to
\Vdual_A$; it takes $v$ to $2\pi i\,q(v,\blank)$. This induces  a morphism $i_q
: \J(\V) \to \Jdual(\V)$ over $\A_g$.

\begin{definition}
The biextension bundle $\Bhat \to \J(\V)$ is the pullback of the bundle
$$
\B(\U) \to \J(\V)\times_{\A_g}\Jdual(\V)
$$
along the map $(\id,i_q): \J(\V) \to \J(\V) \times_{\A_g} \Jdual(\V)$.
The metric $\metric_\Bhat$ is the pullback of the metric $\metric_\B(\U)$.
\end{definition}

\begin{proposition}
\label{integral}
The metric on the biextension bundle $\Bhat$ is constant when restricted to the
zero section, has  translation invariant curvature when restricted to any fiber
$J(V_A)$. It has first Chern class $2\phi = 2\qtilde$. In particular, the class
$2\phi$ is integral.
\end{proposition}

\begin{proof}
The assertions about the metric and its curvature follow immediately from
Proposition~\ref{chern_biext}. We will use Proposition~\ref{chern_biext} to
compute the Chern class.

The class $c_1(\B(\V))$ corresponds to the bilinear form
$$
Q : (V_A\oplus \Vdual_A)^{\otimes 2} \to \C
$$
given by
$$
Q : (u,\phi) \otimes (v,\psi) \mapsto (2\pi i)^{-1}(\phi(v) - \psi(u)).
$$
Pulling this back along the map $(\id,i_q) : V_A \to V_A \oplus \Vdual_A$
we obtain the bilinear form
$$
(u,v) \mapsto
Q\left((u,q(u,\blank)), (v,q(v,\blank)\right) = 2q(u,v).
$$
\end{proof}

\begin{remark}
One can give a moduli interpretation of the points of $\Bhat^\ast$. Points of
$\Bhat_A^\ast$ are biextensions $B$ with weight graded quotients identified with
$\Z$, $V_A$ and $\Z(1)$ that are essentially self dual in the following sense.
If one pulls back the extension
$$
\Hom_\Z(B/\Z(1), \Z(1)) \in \Ext^1_\H(\Vdual, \Z(1))
$$
along the map $i_q: V \to \Vdual$ given by the polarization, then one
obtains $W_{-1}B$.
\end{remark}

The fiber $\Bhat_A^\ast$ of $\Bhat^\ast$ over $[A]\in\A_g$ is a locally
homogeneous space. Let $G_\Z$ be the extension of $V_\Z$ by $\Z$ determined by
the class $2\qtilde \in H^2(V_\Z,\Z)$. So, as a set, $G_\Z = V_\Z \times \Z$.
The multiplication is defined by 
$$
(u,n)\cdot (v,m) = (u + v, m + n + q(u,v)).
$$
Define $G$ to be the complexification of this group. (The underlying set
is $V_\C \times \C$.) Let $F^0G$ by the closed subgroup of $G$ whose underlying
set is $F^0 V\times 0$. Then it is not hard to see that
$$
\Bhat_A^\ast = G_\Z \bs G / F^0G.
$$

In fact, the total space of $\B(\V)^\ast$ and $\Bhat^\ast$ are locally
homogeneous varieties with fundamental groups $\Sp_g(\Z)\ltimes G_{U\, \Z}$
and $\Sp_g(\Z) \ltimes G_\Z$, respectively, as can be seen using
\cite[\S4]{carlson}.

We conclude with a short computation we shall need. Define the commutator
$[a,b]$ of two elements $a$ and $b$ of a group to be $aba^{-1}b^{-1}$.

\begin{lemma}\label{com}
The commutator of the two elements $(u,n)$ and $(v,m)$
of $G_\Z$ is $\bigr(0,2q(u,v)\bigr)$. \qed
\end{lemma}

\section{Theorem~\ref{main} implies Theorem~\ref{asym}}
\label{overview}

In this section, we show that Theorem~\ref{asym} is a direct consequence of
Theorem~\ref{main}. The main ingredient is a fundamental result of Dale Lear
\cite{lear}.

\subsection{Lear's Theorem}

The following is a special case of \cite[(6.2)]{lear} adapted to our present
needs. Suppose that $T$ is a Riemann surface and that $\cC \to T$ is a
semi-stable family of curves over it. This is classified by an orbifold mapping
$f : T \to \Mbar_g$. Denote $f^{-1}(\Mtilde_g)$ by $T'$.

\begin{theorem}
\label{lear}
There exists an integer $N > 0$, which depends only on $g$, such that the
line bundle $f^\ast\B^{\otimes N}$ over $T'$ extends canonically to $T$ as
a holomorphic line bundle with continuous metric. \qed
\end{theorem}

For us, the significance of this result is:

\begin{corollary}
There is an integer $N > 0$ such that the biextension line bundle $\B^{\otimes
N}$ over $\Mtilde_g$ extends naturally to a line bundle on $\Mbar_g$. This
extension is characterized by the requirement that the biextension metric
extend continuously to a metric on the restriction of this line bundle to
to any disk transverse to $\D_0$ at a smooth point.
\end{corollary}

\begin{proof}
Lear's Theorem implies that the integer $N$ and the restriction of the extended
bundle $\Bbar$ to the disk $\disk$ depend only on the local monodromy
representation $\pi_1(\disk^\ast) \to \Sp_g(\Z)\ltimes G_\Z$. This implies that
the extension and the integer are independent of the choice of the disk, which
implies that the extension exists.
\end{proof}

Whether or not the metric on  $\B$ over $\Mtilde_g$ extends to a continuous
metric on $\B$ over $\Mbar_g$ should follow from a good understanding of the
singularities of period mappings of real variations of mixed Hodge structure
in one and several variables. Hopefully the required results will come out of
Pearlstein's work.

We shall regard the extended line bundle as an element of $(\Pic \Mbar_g)
\otimes \Q$. Denote the element that extends $\B$ by $\Bbar$. We have
$$
c_1(\Bbar) =  (8g+4)\lambda + \sum_{h=0}^{[g/2]} r_h \d_h
$$
where each $a_g \in \Q$. We shall compute the coefficients $r_h$ in subsequent
sections.

An immediate consequence of the following result is that Theorem~\ref{main}
implies Theorem~\ref{asym}.

\begin{proposition}
\label{reduction}
With notation as in Theorem~\ref{asym},
$$
\b_g(t) \sim
\begin{cases}
r_0 \log |t| - (4g+2) \log\log\frac{1}{|t|} & \text{ when } h = 0; \cr
r_h \log |t| & \text{ when } h > 0.
\end{cases}
$$
\end{proposition}

\begin{proof}
Let $e$ be a meromorphic section of
$$
\N := \Bbar\otimes\L^{\otimes(-(8g+4))}
$$
that trivializes it over $\M_g$. It is unique up to a non-zero constant. (If
not all the $r_h$ are integers, raise this line bundle to a sufficiently high
power to clear denominators and take $N$th roots later.) Since
$$
c_1(\N) = \sum_{h=0}^{[g/2]} r_h \d_h,
$$
$e$ vanishes to order $r_h$ along $\D_h$. This line bundle is naturally
metrized away from $\D_0$ by the tensor product metric, which we denote by
$\metric_\N$. Since $\metric_\N$ extends smoothly over $\D_h$ when $h>0$, it
follows that if $\disk$ is a small disk with parameter $t$ transverse to a
smooth point of $\D_h$ which it intersects at $t=0$ (as in Theorem~\ref{asym}),
then
$$
\b(t) \sim \log |e(t)|_\N \sim \log|t^{r_h}| = r_h \log|t|
$$
as $t \to 0$.

The situation is slightly more complicated when $h=0$ as the metric on $\L$
does not extend (even continuously) across $\D_0$. To proceed, we compute the
behaviour of the the metric $\metric_\L$ on $\L$ near $\D_0$. This is the
source of the $\log\log(1/|t|)$ term.

Suppose that $\disk$ is a small analytic disk in $\Mbar_g$ with parameter
$t$ that intersects $\D_0$ transversally at a smooth point when $t=0$. Suppose
that $s$ is a trivializing section of the restriction of $\L$ to $\disk$.

\begin{lemma}
\label{asym_L}
As $t \to 0$
$$
\log |s(t)| \sim \frac{1}{2} \log\log \frac{1}{|t|}.
$$
\end{lemma}

\begin{proof}[Proof of Lemma]
Suppose that $C$ is a smooth curve of genus $g$. We use the definition of the
metric on $H^0(C,\Omega^1)$ given in the introduction. Suppose that $a_1,\dots,
a_g,b_1,\dots,b_g$ is a symplectic basis of $H_1(C,\Z)$ and that $w_1,\dots ,
w_g$ is the corresponding normalized basis of $H^0(C,\Omega^1)$. With respect
to these data, the period matrix $\Omega$ of $C$ is given by
$$
\Omega_{jk} = \int_{b_j} w_k = \int_{b_k} w_j,
$$
Thus
$$
w_j = a_j^\ast + \sum_{m} \Omega_{jm} b_m^\ast
$$
where $a_1^\ast,\dots,a_g^\ast,b_1^\ast, \dots, b_g^\ast$ is the dual basis
of $H^1(C,\Z)$. Since it is also symplectic,
$$
(w_j,w_k) = \frac{i}{2} \int_C w_j \wedge \overline{w}_k
= \Im \Omega_{jk}.
$$
Standard linear algebra implies that
$$
\| w_1 \wedge \dots \wedge w_g \|^2 = \det \Im \Omega.
$$

One can construct a trivializing section $s$ over any disk transverse to $\D_0$
at a smooth point as follows: Such a disk corresponds to a stable degeneration
$\{C_t\}_{t\in\disk}$ of smooth genus $g$ curves to an irreducible stable curve
with one node. By \cite[p.~53]{fay}, the normalized period matrix $\Omega(t)$
of $C_t$ is, up to a bounded matrix-valued function of $t$,
$$
\begin{pmatrix}
\Omega_0 & v^T \cr
v & \log t /2 \pi i \cr
\end{pmatrix}
$$
where $\Omega_0$ is the period matrix of the normalization of $C_0$ and $v \in
\C^{g-1}$. It follows that
$$
\det \Im \Omega(t) \sim \text{constant}\cdot \log\frac{1}{|t|}
\text{ as } t \to 0
$$
and that
$$
\| w_1(t) \wedge \dots \wedge w_g (t)\| = \log \det \Im \Omega(t)
\sim \frac{1}{2} \log\log \frac{1}{|t|} \text{ as } t \to 0.
$$
Here $w_1(t), \dots,w_g(t)$ is a normalized basis of the abelian differentials
on the curve $C_t$. Since $w_1(t) \wedge \dots \wedge w_g (t)$ extends to a
local framing of $\L$ on $\Mbar_g$, this establishes the result.
\end{proof}

We can now compute the asymptotics of $\b$ near $\D_0$. Continuing with the
notation from above, we can write the section $e(t)$ of the restriction of $\N$
to the disk $\disk$ in the form
$$
e = t^{r_0} s_B \otimes s_L^{-(8g+4)}
$$
where $s_B$ is a trivializing section of the restriction of $\Bbar$ to $\disk$
and $s_L$ is a trivializing section of the restriction of $\L$ restricted to
$\disk$. Since the metric on $\B$ extends to a continuous metric on the
restriction of $\Bbar$ to $\disk$, Lemma~\ref{asym_L} implies
\begin{align*}
\b(t)
&= r_0\log|t| + \log|s_B(t)|_\B - (8g+4)\log|s(t)|_\L \cr
&\sim r_0\log|t| - (4g+2)\log\log(1/|t|)
\end{align*}
as $t \to 0$.
\end{proof}

\section{A Fundamental Representation}
\label{fund_rep}

The basic idea behind computing the coefficients $r_h$ when $h>0$ is the
following: Suppose that $X$ is a complex manifold and that $D$ is a divisor in
$X$ with irreducible components $D_1,\dots,D_r$. Suppose that $\E$ is a
holomorphic line bundle over $X$ and that $e$ is a meromorphic section of $\E$
that is holomorphic over $X-D$ and trivializes it there. Suppose that
$\gamma_j$ is the boundary of a small holomorphic disk $\disk$ in $X$ that is
transverse to $D_j$. The restriction of the $\C^\ast$ bundle $\E^\ast$ to
$\disk$ is trivial, and thus has fundamental group $\Z$, where the generator
corresponds to a loop in one fiber that encircles the zero section positively.

\begin{lemma}
\label{twist}
The loop $e(\gamma_j)$ equals $n_j$ times the positive generator of
$\pi_1(\E|_\disk^\ast)$ if and only if $e$ has order $n_j$ along $D_j$. \qed
\end{lemma}

In our case we take $\E$ to be the line bundle $\B\otimes
\L^{\otimes(-(8g+4))}$ over $\Mtilde_g$. In this case, we consider the
global representation as the structure of the mapping class group will help
us show that $r_h= -4h(g-h)$.

It is helpful to view the section $e$ as a lift of the normal function $\nu$
as in the following diagram:
$$
\xymatrix{
 & \bigr(\Bhat \otimes \L^{\otimes(-(8g+4))}\bigr)^* \ar[d] \cr
 & \J(\V) \ar[d] \cr
\M_g \ar[r] \ar[ur]^\nu \ar[uur]^{e} & \A_g
}
$$

The framing $e$ therefore induces a representation
$$
\rhohat : \G_g \to
\pi_1\bigr(\bigr(\Bhat \otimes \L^{\otimes(8g+4)}\bigr)^*,*\bigr)
$$
of the mapping class group.

The projection induces a natural homomorphism
$$
\pi_1\bigr(\bigr(\Bhat \otimes \L^{\otimes(8g+4)}\bigr)^*,*\bigr)
\to \pi_1(\A_g,*) \cong \Sp_g(\Z).
$$

\begin{proposition}
The fundamental group of $\bigr(\Bhat \otimes \L^{\otimes(8g+4)}\bigr)^*$
is an extension
$$
1 \to G_\Z \to
\pi_1\bigr(\bigr(\Bhat \otimes \L^{\otimes(8g+4)}\bigr)^*,*\bigr)
\to \Sp_g(\Z) \to 1
$$
Moreover, the composition of the representation $\rhohat$ with the projection
to $\Sp_g(\Z)$ is the natural homomorphism $\G_g \to \Sp_g(\Z)$.
\end{proposition}

\begin{proof}
The kernel is the fundamental group of the pullback of
$\bigr(\Bhat \otimes \L^{\otimes(8g+4)}\bigr)^\ast$ to $\h_g$, the Siegel
upper half space.
Since the pullback of $\L$ to $\h_g$ is trivial, this has the same homotopy
type as the pullback of $\Bhat^\ast$ to $\h_g$. Since $\h_g$ is contractible,
this has the homotopy type of the fiber $\Bhat_A^\ast$ of $\Bhat^\ast$ over 
$[A]$. Since
$$
\Bhat_A^\ast = G_\Z \bs G_\C /F^0G
$$
and since $G_\C/F^0G$ is contractible, the first assertion follows. The second
assertion is clear.
\end{proof}

As an immediate consequence, we see that the restriction of $\rhohat$ to the
Torelli group $T_g$ gives a homomorphism
\begin{equation}\label{tauhat}
\tauhat : T_g \to G_\Z.
\end{equation}
The composition of this homomorphism with the natural homomorphism
$G_\Z \to V_\Z$ is twice the Johnson homomorphism \cite[(6.3)]{hain:comp}.

\section{Topological Preliminaries}

Before computing $r_h$ when $h>0$, we need to clarify our conventions. Although
the final answer is independent of them, signs in intermediate results do. If
$\alpha$ and $\beta$ are paths in a topological space, then $\alpha\beta$ is
the path obtained by following $\alpha$ and then $\beta$. With this convention,
if $X$ is a contractible  space and $\G$ is a group that acts properly
discontinuously and freely on $X$, then for each choice of a base point $x$ of
$X$, there is a natural group isomorphism between $\pi_1(\G\bs X,x)$ and $\G$.
The group element $g \in \G$ corresponds to the homotopy class of the loop in
$\G\bs X$ that is the homotopy class of the projection of a path in $X$ that
goes from $x$ to $gx$.

Now suppose that $\G$ acts on a topological space $S$ on the left. Then $\G$
acts on $X \times S$ on the left. The quotient is a flat $S$ bundle over
$\G \bs X$.

\begin{proposition}
The monodromy of the bundle $\G\bs(X\times X)$ around a loop in
$\G\bs X$ corresponding to $g \in \G$ is $g^{-1} \in \Aut S$.
\end{proposition}

\begin{proof}
The point $(x,s) \in X \times S$ is parallel transported to $(gx,s)$ along the
path corresponding to $g$. But this point is identified with $(x,g^{-1}u)$
under the group action. The result follows.
\end{proof}

Now suppose that $C$ is a smooth projective curve and that $S$ is the
underlying topological space. Since $\M_g$ is the quotient of Teichm\"uller
space on the left by the mapping class group, there is a natural isomorphism
between the mapping class group $\G_g$ of $S$ and $\pi_1(\M_g,[C])$. Applying
the above in the case where $\G = \Diff^+ S$, the group of orientation
preserving diffeomorphisms of $S$, we obtain:

\begin{corollary}\label{monod}
The isotopy class of the geometric monodromy of the universal curve about a
loop in $\M_g$ based at $[C]$ that represents $\phi \in \G_g$ is
$\phi^{-1}$. \qed
\end{corollary}

\section{Coefficient of $\D_h$ when $h>0$}
\label{interior}

Suppose that $g\ge 3$. Fix a stable curve $C_0$ with one node whose moduli
point lies in $\D_h$. There is a stable family of curves $\{C_t\}_{t\in \disk}$
over the disk whose central fiber is $C_0$ and where each other fiber is
smooth. Choose a base point $t_o$ in the punctured disk. We shall use the
oriented surface $S$ underlying $C_{t_o}$ as our reference surface. We can
identify  $\pi_1(\M_g,[C_{t_o}])$ with the mapping class group $\G_g$ of $S$.
The geometric monodromy of this family is a Dehn twist $\s_h:S\to S$ about a
simple closed curve that separates $S$ into two bounded surfaces, one of genus
$h$, the other of genus $g-h$.

Note that $\s_h$ lies in the Torelli group $T_g$ and, moreover, that it lies
in the kernel of the Johnson homomorphism $T_g \to V_\Z$. Consequently, its
image under the homomorphism $\tauhat : T_g \to G_\Z$, constructed in
Section~\ref{fund_rep}, lies in the center $\Z$ of $G_\Z$.

\begin{proposition}
\label{coeff}
If $g\ge 3$, then $\tauhat(\s_h) = 4h(g-h)$.
\end{proposition}

This computation is equivalent to (up to a scaling constant) to a result of
Morita \cite[Thm.~5.7]{morita:casson} and will be proved below. Combining it
with Proposition~\ref{monod} completes the computation of $r_h$ when $h>0$.

\begin{corollary}
If $h>0$, then $r_h = -4h(g-h)$. \qed
\end{corollary}

\begin{proof}[Proof of Proposition~\ref{coeff}]
Let $W$ be a genus $h$ subsurface of the reference surface $S$ such
that $\s_h$ is the Dehn twist about the boundary of $W$. Choose a base point
$\ast \in \del W$. The boundary of $W$ (with the induced orientation)
determines an element $c$ in $\pi_1(W,\ast)$. One can choose generators
$a_1,\dots, a_h, b_1,\dots, b_h$ of $\pi_1(W,\ast)$ whose homology classes
comprise a symplectic basis of $H_1(W,\del W)$ and such that
$$
\prod_{j=1}^h [a_j,b_j] = c^{-1}.
$$

We can regard $W$ as the real oriented blow up of the compact surface $W'$
obtained by identifying $\del W$ to a point $P$. We shall use $P$ as a 
base point of $W'$. Let $F$ be the unit tangent bundle of $W'$. The base point
$\ast$ of $W$ corresponds to a point in the fiber of $F$ over $P$
that we shall also denote by $\ast$ and shall use as a base point for $F$.

Choose a vector field $\xi$ on $W'$ whose only zero is at $P\in W'$. This
induces a map $\xi : (W,\ast) \to (F,\ast)$ and therefore a homomorphism
$$
\xi_\ast : \pi_1(W,\ast) \to \pi_1(F,\ast).
$$
Denote by $f$ the element of $\pi_1(F,\ast)$ that corresponds to going once
around the fiber over $P$ in the positive (i.e., counterclockwise) direction.
It follows from the Hopf Index Theorem that $\xi_\ast (c) = f^{2h-2}$.
Consequently
$$
\xi_\ast\left(\prod_{j=1}^h [a_j,b_j]\right) = f^{2-2h}.
$$

Denote by $\Aut_0 \pi_1(W,*)$ the group of automorphisms $\phi$ of $\pi_1(W,*)$
such that $\phi(c) = c$. There is a natural homomorphism
$$
\chi : \pi_1(F,*) \to \Aut_0 \pi_1(W,*)
$$
that covers the natural left action of $\pi_1(W')$ on itself by conjugation.
Note that $\s_h$ can be viewed as an element of $\Aut_0 \pi_1(W,\ast)$. In
fact,
$$
\s_h : x \mapsto c x c^{-1} = \chi(f^{-1})(x).
$$
It follows that
$$
\s_h^{2h-2} = \chi(f^{2-2h}) =
\chi\circ\xi_\ast\left(\prod_{j=1}^h [a_j,b_j]\right) \in \Aut_0 \pi_1(W,*).
$$

Observe that the image of $\pi_1(W,*)$ in $\Aut_0 \pi_1(W,\ast)$ is naturally a
subgroup of $\G_g$ as each automorphism coming from $\pi_1(W, *)$ is induced by
a diffeomorphism of $W$ that fixes $\del W$ pointwise. In fact, the image lies
in the Torelli group. Let
$$
\phi : \pi_1(W,*) \to T_g
$$
be the composite. Then, since $G_\Z$ is a two-step nilpotent group
$$
\tauhat(\s_h) =
\frac{1}{2h-2}\,\tauhat\left(\prod_{j=1}^h [\phi(a_j),\phi(b_j)]\right)
= \frac{1}{2h-2}\,
\sum_{j=1}^h \,[\tauhat\circ\phi(a_j),\tauhat\circ\phi(b_j)] \in G_\Z
$$
and the value of this expression depends only on the images
of the $\phi(a_j)$ and $\phi(b_j)$ under the Johnson homomorphism
$T_g \to V$ (cf.\ Lemma~\ref{com}).

In order to compute the image of $\s_h$ in $G_\Z$, we shall need the following
result.
Denote the homology class of $x\in \pi_1(W,*)$ by $[x]$.
Let $w' = \sum_{j=1}^h [a_j] \wedge [b_j]$. We view this as an element of 
$\Lambda^2 H_1(W)$.

\begin{lemma}
The image of $x \in \pi_1(W,*)$ under the composite
$$
\begin{CD}
\pi_1(W,*) @>\phi>> T_g @>{\mathrm{Johnson}}>> V,
\end{CD}
$$
where $\phi$ is the composite
$$
\begin{CD}
\pi_1(W,\ast) @>>> \pi_1(C,\ast) @>{\mathrm{inner}}>> T_g,
\end{CD}
$$
is the image of $[x] \wedge w' \in \Lambda^3 H_1(W)$ in $V = \Lambda^3 H_1(S)/
H_1(S)$ under the map induced by the natural inclusion $H_1(W) \hookrightarrow 
H_1(S)$.
\end{lemma}

\begin{proof}
Denote the surface obtained from $W$ by collapsing the boundary to a point by
$\Wbar$. Observe that the quotient mapping induces an isomorphism on $H_1$. Set
$$
\G_{W,\del W} = \pi_0 \Diff^+ (W,\del W) \text{ and }
\G_{\Wbar,\ast} = \pi_0\Diff^+ (\Wbar,\ast),
$$
where $\ast$ is the image of $\del W$ in $\Wbar$, be the mapping class groups
of $(W,\del W)$ and $(\Wbar,\ast)$, respectively. Denote the associated Torelli
groups by $T_{W,\del W}$ and $T_{\Wbar,\ast}$. In Harer's notation, these are
isomorphic to $T_{h,1}$ and $T_h^1$, respectively. There are Johnson
homomorphisms
\begin{equation}
\label{double}
T_{W,\del W} \to \Lambda^3 H_1(W) \text{ and }
T_{\Wbar,\ast} \to \Lambda^3 H_1(\Wbar)
\end{equation}
which are compatible in the sense that the first is the composition of the
second with the natural quotient mapping $T_{W,\del W} \to T_{\Wbar,\ast}$ once
one identifies $H_1(W)$ with $H_1(\Wbar)$. The homomorphisms (\ref{double})
both induce isomorphisms on $H_1$ mod torsion.

The first step is to prove that the composite $$ H_1(\pi_1(\Wbar,\ast)) \to
H_1(T_{\Wbar,\ast}) \to \Lambda^3 H_1(\Wbar) $$ takes $[x]$ to $[x]\wedge \w'$.
This follows directly from the geometric definition of the Johnson homomorphism
that is given in \cite{johnson:survey} and which is explained in detail in
\cite[\S3]{hain:normal}. Indeed, the mapping torus $M \to S^1$ that corresponds
to conjugation by $x$ is trivial as the diffeomorphism $\phi$ of $C$ which
induces it is isotopic to the identity when the base point is ignored. The
section of base points takes $t \in S^1$, which we regard as $[0,1]/0\sim 1$,
to $\gamma(1-t)$, where $\gamma$ is any loop representing $x$. The map from $M$
to $\Jac C$ then takes $(u,t) \in C\times S^1$ to $u - \gamma(1-t)$. This has
homology class $[x]\wedge \w'$. 

By the compatibility of the Johnson homomorphisms for $(W,\del W)$ and
$(\Wbar,\ast)$, it follows that
$$
H_1(T_{W,\del W}) \to \Lambda^3 H_1(W)
$$
takes $[x]$ to $[x]\wedge \w'$.

The result now follows as the inclusion of $W$ into $S$ induces a homomorphism
$T_{W,\del W} \to T_g$ such that the diagram
$$
\begin{CD}
T_{W,\del W} @>>> T_g \cr
@VVV @VVV \cr
\Lambda^3 H_1(W) @>>> V
\end{CD}
$$
commutes, where the bottom arrow is induced by the induced map $H_1(W) \to
H_1(S)$.
\end{proof}

Since the composite of $\tauhat: T_g \to G_\Z$ with the natural
homomorphism $G_\Z \to V$ is twice the Johnson homomorphism, it follows
from Lemma~\ref{com} that
$$
\tauhat(\s_h) = \frac{8}{2h-2}\,
\sum_{j=1}^h \, q\bigr([a_j]\wedge w',[b_j]\wedge w'\bigr) \in \Z.
$$
We will abuse notation and not distinguish between 
an element of $\pi_1(W,*)$ and its homology class. We can extend the
symplectic basis $a_1,\dots, r_h, b_1, \dots, b_h$ of $H_1(W)$ to a 
symplectic basis $a_1,\dots, a_g, b_1, \dots, b_g$ of $H_1(S)$. Set
$$
w = \sum_{j=1}^g a_j \wedge b_j
$$
and $w'' = w - w'$.
Denote $H_1(W)$ by $H'$ and its orthogonal complement in $H_1(S)$
by $H''$. Recall from Section~\ref{lin_alg} that if $u,v \in V$, then 
$q(u,v) = \bil {j(u)} {j(v)}/(g-1)$. A short computation shows that if
$x \in H'$, then
$$
j(x\wedge w') = (g-1)\, x \wedge w' - (h-1)\, x \wedge w 
= (g-h)\, x \wedge w' - (h-1)\, x \wedge w''.
$$
Since $H'\wedge w''$ is orthogonal to $H''\wedge w'$ in $\Lambda^3 H_1(S)$,
we have
\begin{equation*}
\begin{split}
&\sum_{j=1}^h \,q(a_j\wedge w',b_j\wedge w') \cr
=&\frac{1}{g-1}\,\sum_{j=1}^h
\left( (g-h)^2\bil {a_j\wedge w'} {b_j\wedge w'} +
(h-1)^2 \bil {a_j\wedge w''} {b_j\wedge w''}\right) \cr
=& \frac{1}{g-1}\,\sum_{j=1}^h \bigr((g-h)^2(h-1) + (h-1)^2(g-h)\bigr) \cr
=& \frac{1}{g-1}\,h(h-1)(g-h)(g-1)\cr
=& h(h-1)(g-h) \cr
\end{split}
\end{equation*}
Thus $\tauhat(\s_h) = 8h(g-h)(h-1)/(2h-2) = 4h(g-h)$.
\end{proof}

\begin{corollary}
\label{compact_part}
The first Chern class of the biextension bundle $\B$ over $\Mtilde_g$ is
$$
c_1(\B) = (8g+4)\lambda - 4\sum_{h=1}^{[g/2]} h(g-h) \d_h. \qed
$$
\end{corollary}

\section{The Coefficient of $\D_0$}

Our final task is to determine the coefficient of $\d_0$. In principle we could
do this using the representation defined in Section~\ref{fund_rep}. However,
this appears to be quite complicated, and it seems easier to do it directly by
restriction to the hyperelliptic locus, where we know that $\Bbar$ is trivial,
and using a known linear equivalence to show $r_0 = -g$.

The set of points of $\Mbar_g$ corresponding to stable hyperelliptic curves
forms a subvariety $\Hbar_g$. When $h>0$, the intersection $\Delta_j \cap
\Hbar_g$ is an irreducible divisor. We shall denote its class by $\d_j$. When
$h=0$,
$$
\D_0 \cap \Hbar_g = \Xi_0 \cup \bigcup_{h=1}^{[(g-1)/2]} \Xi_h
$$
where $\Xi_h$, $1\le h \le (g-1)/2$ is the closure of the locus in $\Hbar_g$
consisting of the stable hyperelliptic curves that are the union of two smooth
projective curves, one of genus $h$ and another of genus $g-h-1$, which
intersect transversally in two points. The divisor $\Xi_0$ consists of the
closure of the locus of irreducible hyperelliptic curves whose normalization is
of genus $g-1$. Each $\Xi_h$ is a divisor in $\Hbar_g$.  Denote the class of
$\Xi_h$ in $\Pic \Hbar_g$ by $\xi_h$.

Cornalba and Harris \cite[Prop.~4.7]{corn-har} prove that the classes $\xi_h$,
$0\le h \le (g-1)/2$ and $\d_h, 1\le h \le g/2$ are linearly independent in
$\Pic \Hbar_g$ and that in $\Pic\Hbar_g$, one has
$$
(8g+4)\lambda = g\xi_0 + \sum_{j=1}^{[(g-1)/2]}(j+1)(g-j) \xi_j
+ 4\sum_{h=1}^{[g/2]}h(g-h) \d_h.
$$

Proposition~\ref{hyperelliptic} implies that $\nu^\ast \phi$ vanishes on
$\Hbar_g$. Thus, we have
$$
(8g+4)\lambda =
-r_0\xi_0 + 4\sum_{h=1}^{[g/2]} h(g-h)\d_h + \sum_{j\ge 1} c_j \xi_j.
$$
in $H^2(\Hbar_g,\Q)$ and thus in $\Pic \Hbar_g$. Since, when $j>0$, $\Xi_j$ is
contained in the codimension 2 stratum of $\Mbar_g$, it follows that $r_0 =
-g$, as claimed. This completes the proof of Theorem~\ref{main}.

\end{document}